\documentclass[11pt, letterpaper]{amsart}
\usepackage{fullpage}
\usepackage{microtype}
\usepackage[dvipsnames]{xcolor}
\usepackage[croatian,english]{babel}
\usepackage{colonequals}
\usepackage{appendix}
\usepackage{amsmath,amssymb}
\usepackage{varwidth}
\usepackage{versions}
\usepackage{tikz-cd}
\usepackage{mathrsfs}
\usepackage{enumitem}
\usepackage[colorlinks=true,hyperindex,pagebackref,linktocpage=true]{hyperref}
\usepackage{cleveref}
\excludeversion{misc}

\hypersetup{
	colorlinks,
	linkcolor={teal!60!blue},
	citecolor={Goldenrod!60!Maroon},
	urlcolor={purple!80!black}
}

\newtheorem{thm}{Theorem}[section]
\newtheorem{prop}[thm]{Proposition}
\newtheorem{lem}[thm]{Lemma}

\newtheorem{defn}[thm]{Definition}
\newtheorem{cor}[thm]{Corollary}
\newtheorem{rem}[thm]{Remark}
\newtheoremstyle{named}{}{}{\itshape}{}{\bfseries}{.}{.5em}{#1 \thmnote{#3}}
\theoremstyle{named}

\Crefname{lem}{lemma}{lemmas}
\Crefname{lem}{Lemma}{Lemmas}
\Crefname{thm}{theorem}{theorems}
\Crefname{thm}{Theorem}{Theorems}
\Crefname{defn}{definition}{definitions}
\Crefname{defn}{Definition}{Definitions}
\Crefname{prop}{proposition}{propositions}
\Crefname{prop}{Proposition}{Propositions}
\Crefname{cor}{corollary}{corollaries}
\Crefname{cor}{Corollary}{Corollaries}
\Crefname{conj}{conjecture}{conjectures}
\Crefname{conj}{Conjecture}{Conjectures}

\newcommand{\brF}{\breve{F}}
\newcommand{\brG}{\breve{G}}

\newcommand{\calO}{{\mathcal{O}}}
\newcommand{\calD}{{\mathcal{D}}}
\newcommand{\calC}{{\mathcal{C}}}

\newcommand{\tri}{\Delta}
\newcommand{\q}{\mathbf{q}}
\newcommand{\p}{\mathbf{p}}
\newcommand{\sigen}{\langle\sigma\rangle}
\newcommand{\bgmuindec}{B(G,\mu)_\textnormal{indec}}
\newcommand{\ssigma}{\mathbb{S}/\!\langle \sigma\rangle}
\newcommand{\renv}[1]{A(G,\mu)_{\ge #1}}

\DeclareMathOperator{\area}{A}

\DeclareMathOperator{\indec}{indec}

\DeclareMathOperator{\fin}{fin}

\DeclareMathOperator{\height}{ht}
    \title{Hodge-Newton indecomposability and\\ a combinatorial identity}
	
\begin{document}

	\author[D.G. Lim]{Dong Gyu Lim}
	
	\email{limathonggyu@gmail.com}
	
	\begin{abstract}
	We present a simple alternative viewpoint on Hodge-Newton indecomposability, illustrating its explanatory value through a uniform proof of a combinatorial identity arising from affine Deligne-Lusztig varieties with finite Coxeter part.
	\end{abstract}
	\date{\today}
	
	\maketitle

\section{Introduction}

\subsection{Motivation}

Affine Deligne-Lusztig varieties were first introduced by Rapoport in \cite{Ra05} to study mod $p$ reduction of Shimura varieties. Compared with the hyperspeical level, affine Deligne-Lusztig varieties with Iwahori level structure exhibit more complicated geometry: properties such as nonemptiness, dimension, and equidimensionality are governed by subtle combinatorial conditions. In \cite{HNY}, the authors study a certain class of affine Deligne-Lusztig varieties at Iwahori level and show that their geometric structure is particularly simple. One of the challenges in their analysis is the proof of a combinatorial identity arising from the study of class polynomials; the simplest instance of this identity is the following.

\begin{thm}[{\cite[1C]{HNY}}]\label{motiv}
Fix natural numbers $i< n$, then
\begin{align*}
    \sum_{\substack{k\ge1,1>\frac{a_1}{b_1}>\cdots>\frac{a_k}{b_k}>0\\a_1+\cdots+ a_k=i,b_1+\cdots+b_k=n}} (\q-1)^{k-1}\q^{1-k+\frac{\sum_{1\le l_1<l_2\le k}(a_{l_1}b_{l_2}-a_{l_2}b_{l_1}) +\sum_{1\le l\le k}\gcd(a_l,b_l)}{2}}=\q^{\frac{i(n-i)-n}{2}+1}.
\end{align*}
\end{thm}

The authors ask whether this formula admits a purely combinatorial proof. The argument in \textit{loc.cit.} relies mainly on the Chen-Zhu conjecture for affine Deligne-Lusztig varieties at the hyperspecial level, which appears to be a rather heavy tool for this purpose. Moreover, this identity is only a special case -- namely, the (split) root system $A_{n-1}$ with the $i$-th fundamental coweight -- of a more general formula of the following type (see \Cref{spadesuit} for the precise statement).
\begin{equation}\tag{$\star$}\label{generalform}
    \sum_{[b]\in \bgmuindec} (\q-1)^?\q^{-??}=1.
\end{equation}

This general formula is established using a combination of additional tools, including extensive case-by-case computations and computer assistance. While fully rigorous, this approach does not immediately highlight the combinatorial structure, particularly the role of the index set and its relation to the exponents, is not immediately apparent. In this paper, we offer a new perspective on $\bgmuindec$ that provides a natural interpretation of the formula and clarifies the exponents.

\subsection{Convex Hulls}\label{convexhull}

We briefly outline the proof of \Cref{motiv} to illustrate the main idea: \textit{the index set in the formula can be interpreted as a set of certain convex hulls}.

Let $\calD$ denote the index set. For a natural number $k$, let $\calD_k$ be the subset consisting of elements with $k$ tuples. Replacing $b_i$'s by $b_i-a_i$, we obtain another index set $\calC_k$, consisting of tuples $(x_l,y_l)\in\mathbb{Z}_{>0}^2$ for $l\le k$ such that \[x_1 + \cdots + x_k = i\text{, }\quad y_1+\cdots+y_k = n-i\text{,\quad and }\quad\frac{y_1}{x_1} < \cdots <\frac{y_k}{x_k}.\] In the coordinate plane $\mathbb{R}^2$, the set $\calC_k$ can be identified with the set of all convex\footnote{To avoid confusion, this is convex upward.} broken lines connecting $O\colonequals(0,0)$ and $Y\colonequals(i,n-i)$ with $k$ segments and lying strictly above the broken line $\overline{OXY}$, where $X\colonequals(i,0)$ (see \Cref{btob-a}). After the substitution $\p \longleftrightarrow\frac{\q-1}{\q}$ and an appropriate interpretation of the exponents, the problem reduces to proving the following proposition, where $\calC = \bigcup_{k\ge1}\calC_k$ (see \Cref{pick}.)

\begin{prop}\label{essentialidea}
The following identity holds:
\[\sum_{C\in \mathcal{C}}(1-\p)^{u(C)} \p^{b(C)}=1,\]
where $u(C)$ denotes the number of lattice points lying strictly under $C$, and $b(C)$ the number of break points of $C$, both lying in the interior of the $\triangle OXY$.
\end{prop}

\begin{proof}
    We may assume $0<\p<1$. Consider the following random process.
    
    For each lattice point in the interior of $\triangle OXY$, select it independantly with probability $p$. Let $P$ be the convex hull of all selected points together with $O$ and $Y$. Removing the segment $\overline{OY}$, the boundary of $P$ determines a broken line $C\in\calC$ connecting $O$ and $Y$. Moreover, a given $C\in\calC$ occurs exactly when

    a) the lattice points strictly under $C$ are not selected, and

    b) the break points of $C$ are selected.

Hence, the probability of this event is $(1-\p)^{u(C)}\p^{b(C)}$. Summing over all $C\in\calC$, the total probability is $1$.
\end{proof}

The proof of the general version (see \Cref{maincor}) is essentially the same. The index set of \eqref{generalform} is interpreted as a collection of convex hulls, and the exponents count lattice points under (or over) the hull and on the hull.

\subsection{Notations}

Let $F$ be a nonarchimedean local field and $G$ a connected reductive group over $F$. Let $\brF$ be the completion of the maximal unramified extension of $F$, write $\brG = G(\brF)$. The Frobenius map on the residue field of $\brF$ lifts to an automorphism $\sigma$ of $\brF$, which we extend to $\brG$. Let ${\Gamma_0}$ be the absolute Galois group of $\brF$.

Let $S$ be a maximal $\brF$-split torus of $G$ over $F$, and $T$ its centralizer. Since $G$ becomes quasi-split over $\brF$, $T$ is a maximal torus. Set $V \colonequals X_*(T)_{\Gamma_0,\mathbb{Q}}$. Fix a base alcove in the apartment of $S$, which determines simple coroots $\alpha_i^\vee$ and the set $\mathbb{S}$ of simple reflections. The automorphism $\sigma$ acts on $V$ and hence on $\mathbb{S}$. For $i\in\mathbb{S}$, let $\alpha_i$, $\varpi_i$, and $\varpi_i^\vee$ denote the corresponding root, fundamental weight, and fundamental coweight. Let $\calO_i$ as the $\sigma$-orbit of $i$, and set $\varpi_{\calO_i}\colonequals\sum_{j\in\calO_i}\varpi_{j}$.

We assume that $G$ is unramified and adjoint, and that $\mu$ is not central on any $\sigma$-orbit of connected components of the $\brF$-root system of $G$. These technical assumptions mainly exclude trivial cases; see the remark below for justification. Under these assumptions, the general formula \eqref{generalform} reduces to the following identity (\cite[6A]{HNY} $\spadesuit'$), which will be proved in \Cref{spadesuit}.
\[\sum_{[b]\in \bgmuindec}(\q-1)^{\#(\ssigma)-\# (I(\nu_b)/\!\sigen)}\q^{-\sum_{i\in \ssigma}\lceil\langle\mu-\nu_b,\varpi_{\calO_i}\rangle\rceil+\# (I(\nu_b)/\!\sigen)} =1.\]

\begin{rem}
We assume that $G$ is unramified so that the relative root system over $\brF$ is reduced, and that $G$ is adjoint, both in order to apply \cite[Lemma 3.5]{HeNie}. This lemma states that \Cref{defbgmu} is the bijective image of the standard $B(G,\mu)$ under the Newton map.

If there exists a $\sigma$-orbit $\Phi$ of connected components of the root system where $\mu$ is central, the ``$\Phi$-part'' in each exponent of \eqref{generalform} will be simply zero. Moreover, the ``$\Phi$-part'' of $\bgmuindec$ is a singleton. Therefore, in this case, \eqref{generalform} reduces to the formula with $B(G',\mu')_\textnormal{indec}$ where $G'\le G$ and $\mu'$ is not central on any $\sigma$-orbit of connected components of the root system of $G'$.
\end{rem}



\section{$ B(G,\mu)_\textnormal{indec}$ as the set of convex hulls}\label{indeconvex}

\subsection{$B(G,\mu)$ and lattice points}
For our purposes, the following definition of $B(G,\mu)$ will suffice. 

\begin{defn}\label{defbgmu}
    Let $\mu \in (V^+)^\sigma$. For $v\in V$, define $I(v)\colonequals \{i\in\mathbb{S}: \langle v,\alpha_i\rangle=0\}$ and \[B(G,\mu)\colonequals \{v\in (V^+)^\sigma:\langle \mu-v,\varpi_{\calO_i}\rangle\in\mathbb{Z}_{\ge0}~\forall i\in \mathbb{S}\setminus\! I(v)\}.\]
\end{defn}

\begin{defn}
Let $\mu\in (V^+)^\sigma$ and $J\subset \mathbb{S}$ be $\sigma$-stable. For $v\in B(G,\mu)$, we say the pair $(\mu,v)$ is HN-decomposable with respect to $J$ if $I(v)\subset J$ and $\mu-v \in \sum_{j\in J}\mathbb{Q}_{\ge0}\alpha_j^\vee$.

If no such $\sigma$-stable $J\subsetneq\mathbb{S}$ exists, the pair $(\mu,v)$ is called HN-indecomposable and the set of all such $v\in B(G,\mu)$ is denoted by $\bgmuindec$.
\end{defn}

As discussed in \Cref{convexhull}, our goal is to understand $\bgmuindec$ as a set arising from certain convex hulls. We begin by introducing the relevant lattice points.

\begin{defn}[Set of ``lattice'' points]
    In $\ssigma \times \mathbb{Q}_{>0}$, we define the set $L$ by
    \[L =\left\{(\calO_i,\langle\mu,\varpi_{\calO_i}\rangle-m):\calO_i\in\ssigma, m\in\mathbb{Z}, 1\le m<\langle\mu,\varpi_{\calO_i}\rangle\right\}.\]
\end{defn}

Note that if $\langle\mu, \varpi_{\calO_{j}}\rangle < 1$ for some $\calO_{j}\in \ssigma$, then there is no element of $L$ of the form $(\calO_{j}, -)$.

\begin{defn}
    Let $L_0\subset L$. For each $\calO_i\in\ssigma$, define $\height_{L_0}(\calO_i) \colonequals \max \{y: (\calO_i, y)\in L_0\}$ whenever this set is nonempty; otherwise we set $\height_{L_0}(\calO_i) \colonequals 0$.
\end{defn}

\begin{defn}
For a subset $L_0\subset L$, define\[B(G,\mu)_{\ge L_0} \colonequals \{v \in B(G,\mu): \langle v, \varpi_{\calO_i}\rangle \ge \height_{L_0}(\calO_i) \text{ for all }\calO_i\in\ssigma\}.\]
\end{defn}

The central observation of this paper is the following theorem.

\begin{thm}\label{indecsup}
    The set $\bgmuindec$ consists precisely of the smallest elements of $B(G,\mu)_{\ge L_0}$, as $L_0$ ranges over all subsets of $L$.
\end{thm}

\begin{proof}
    Combine \Cref{equivsmall} and \Cref{indecisminofL0} below.
\end{proof}

\begin{prop}\label{equivsmall}
    Let $v\in B(G,\mu)_{\ge L_0}$. The following are equivalent:
    \begin{enumerate}
        \item[1)] $v$ is the smallest element of $B(G,\mu)_{\ge L_0}$.
        \item[2)] $v$ is minimal in $B(G,\mu)_{\ge L_0}$.
        \item[3)] For all $i \in \mathbb{S}\setminus\! I(v)$, we have $\langle v,\varpi_{\calO_i}\rangle= \height_{L_0}(\calO_i)$.
    \end{enumerate}
    
    The set $B(G,\mu)_{\ge L_0}$ is always nonempty, since it contains $\mu$. Hence, by taking any minimal element, we obtain the smallest element, which we denote by $\min B(G,\mu)_{\ge L_0}$. Moreover, we claim that $\min B(G,\mu)_{\ge L_0}\in\bgmuindec$.
\end{prop}

\Cref{equivsmall} closely resembles Theorem 6.5 in \cite{Chai}, but requires additional arguments, including integer conditions and a bounded set instead of the full chamber. We also claim new properties, such as those concerning $\bgmuindec$. Moreover, our approach differs considerably. The full proof is given in \Cref{proofprop}.

\begin{lem}\label{indecisminofL0}
Let $v'\in\bgmuindec$ and define $L_{v'} \colonequals \{(\calO_i, \langle v',\varpi_{\calO_i}\rangle): i \in \mathbb{S}\setminus\! I(v')\}$. Then\[\min B(G,\mu)_{\ge L_{v'}} = v'.\]
\end{lem}

\begin{proof}
By definition of $L_{v'}$, we have $\langle v', \varpi_{\calO_i} \rangle = \height_{L_{v'}}(\calO_i)$ for all $i\in \mathbb{S}\setminus\! I(v')$. Hence condition 3) of \Cref{equivsmall} is satisfied.
\end{proof}

\subsection{Proof of the general formula \eqref{generalform}}
As a corollary, we naturally obtain the following identity.

\begin{cor}\label{maincor}
Let $\mu\in V^+$ be non-central on every $\sigma$-orbit of the connected components of the root system. Then, as polynomials in a variable $\p$,\[\sum_{v\in B(G,\mu)_{\textnormal{indec}}}(1-\p)^{\sum_{\calO_i\in\ssigma}\lceil\langle \mu-v,\varpi_{\calO_i}\rangle\rceil-\#(\ssigma)} \p^{\#(\ssigma)-\#(I(v)/\!\sigen)}=1.\]
\end{cor}

\begin{proof}

    We may assume $0<\p<1$ and consider the following process. For each lattice point of $L$, select the point independently with probability $\p$. Let $L_0$ be the set of all selected points. Take the smallest element of $B(G,\mu)_{\ge L_0}$, which exists and lies in $\bgmuindec$ by \Cref{equivsmall}. Moreover, it equals $v$ exactly when, by condition 3) of \Cref{equivsmall}, 

    a) $\height_{L_0}(\calO_i) \le \langle v, \varpi_{\calO_i}\rangle$ for all $i\in \mathbb{S}$ and,

    b) $\height_{L_0}(\calO_i) = \langle v,\varpi_{\calO_i}\rangle$ for all $i \in \mathbb{S}\setminus\! I(v)$.

    The condition a) holds when all lattice points $(\calO_i, y)\in L$ such that $y > \langle v,\varpi_{\calO_i}\rangle$ are not selected. We claim that $\langle \mu - v,\varpi_{\calO_i}\rangle >0$ for all $i\in\mathbb{S}$, so that the number of such lattice points is $\lceil\langle \mu-v,\varpi_{\calO_i} \rangle\rceil - 1$ for all $i\in \mathbb{S}$. Suppose $\langle \mu-v,\varpi_{\calO_k}\rangle = 0$ for some $k\in\mathbb{S}$. Since $v\in \bgmuindec$, this implies $k\not\in I(v)$. Then by condition 3) of \Cref{equivsmall} we have $\langle v,\varpi_{\calO_k}\rangle = \height_{L_0}(\calO_k) < \langle \mu,\varpi_{\calO_k}\rangle$, a contradiction.

    The condition b) holds when the lattice points $(\calO_i,\langle v,\varpi_{\calO_i}\rangle)$ are selected for $i\in \mathbb{S}\setminus\! I(v)$. Note that $\langle v,\varpi_{\calO_i}\rangle < \langle\mu,\varpi_{\calO_i}\rangle$ since $v\in \bgmuindec$, and $\langle \mu,\varpi_{\calO_i}\rangle - \langle v,\varpi_{\calO_i}\rangle \in \mathbb{Z}$ because $v\in B(G,\mu)$. Hence $(\calO_i,\langle v,\varpi_{\calO_i}\rangle)\in L$.

    Therefore the probability of this event is
    \begin{align*}    
    &(1-\p)^{\sum_{\calO_i\in\ssigma}(\lceil\langle \mu-v,\varpi_{\calO_i} \rangle\rceil - 1)}\p^{\sum_{\calO_i\in(\mathbb{S}\setminus\!I(v))/\!\sigen}1} \\=& (1-\p)^{\sum_{\calO_i\in\ssigma}\lceil\langle \mu-v,\varpi_{\calO_i} \rangle\rceil - \#(\ssigma)} \p^{\#(\ssigma)-\#(I(v)/\!\sigen)}.
    \end{align*}

    By \Cref{indecsup}, every element of $\bgmuindec$ can arise from this process. Hence the sum of the above probabilities over $\bgmuindec$ is $1$.
\end{proof}

\begin{cor}[{\cite[6A]{HNY} $\spadesuit'$}]\label{spadesuit}
For $\mu\in (V^+)^\sigma$ non-central on every $\sigma$-orbit of the connected components, we have
    \[\sum_{v\in B(G,\mu)_{\indec}}(\q-1)^{\#(\ssigma)-\# (I(v)/\!\sigen)}\q^{-\sum_{i\in \ssigma}\lceil\langle\mu-\nu_b,\varpi_{\calO_i}\rangle\rceil+\# (I(\nu_b)/\!\sigen)} =1.\]
\end{cor}
\begin{proof}
    Apply $\p=\frac{\q-1}{\q}$ to \Cref{maincor}.
\end{proof}

\subsection{Proof of \Cref{equivsmall}}\label{proofprop}

For $J\subset \mathbb{S}$, denote by $\partial J$ the set of vertices of $\mathbb{S}\setminus\!J$ adjacent to $J$.

\begin{lem}\label{subdiag}
    Let $J\subset\mathbb{S}$ be connected and $j$ be a vertex of $J$. Then, there exist $c_i\in\mathbb{Q}_{>0}$ for each $i\in J$ and $d_k\in \mathbb{Q}_{>0}$ for each $k\in \partial J$ such that\[\varpi_j^\vee=\sum_{i\in J}c_i\alpha_i^\vee+\sum_{k\in\partial J} d_k\varpi_k^\vee.\]
\end{lem}

\begin{proof}
    Noting that $J$ is a Dynkin diagram, we can find the usual coroot coefficients $c_i$'s of the fundamental coweight corresponding to $j$. They are all positive (cf. \cite[Lemma 2.18]{DGL23}). Calculating $d_\ell\colonequals\langle\varpi_j^\vee-\sum_{i\in J}c_i\alpha_i^\vee, \alpha_\ell\rangle$ for each $\ell \in \mathbb{S}\setminus\!J$, we can find the $\sum_k d_k\varpi_k^\vee$ part. It is nonzero only when $\ell\in\partial J$. In that case, it easily follows that $d_\ell$ is positive.
\end{proof}

The following lemma makes us perturb $v\in (V^+)^\sigma$ such that $v \le \mu$ and $\langle v,\varpi_{\calO_k}\rangle \neq \height_{L_0}(\calO_k)$ for some $\calO_k\in\ssigma$ into a more desirable one. For simplicity, we temporarily define\[\renv{L_0} \colonequals\{v\in (V^+)^\sigma: v\le \mu\text{ and } \langle v,\varpi_{\calO_i}\rangle\ge \height_{L_0}(\calO_i)~\forall \calO_i\in\ssigma\}.\]Obviously, $B(G,\mu) \cap \renv{L_0} = B(G,\mu)_{\ge L_0}$ for any $L_0\subset L$.

\begin{lem}\label{goingdown}
Let $L_0\subset L$ and $v\in \renv{L_0}$. Suppose that there exists $k\in \mathbb{S}\setminus\! I(v)$ such that $\langle v,\varpi_{\calO_{k}}\rangle\neq \height_{L_0}(\calO_k)$. Then, one can find $v'\in \renv{L_0}$ such that either \eqref{inS} or \eqref{breakpts} holds.\begin{equation}\tag{S}\label{inS}
\#\{i\in\mathbb{S}\!:\langle v',\varpi_{\calO_i}\rangle=\height_{L_0}(\calO_i)\}>\#\{i\in\mathbb{S}\!:\langle v,\varpi_{\calO_i}\rangle=\height_{L_0}(\calO_i)\}\end{equation}\begin{equation}\tag{I\textsuperscript{c}}\label{breakpts}
     \#\{i\in\mathbb{S}\setminus\! I(v')\!:\langle v',\varpi_{\calO_i}\rangle\neq \height_{L_0}(\calO_i)\}<\#\{i\in\mathbb{S}\setminus\! I(v)\!:\langle v,\varpi_{\calO_i}\rangle\neq \height_{L_0}(\calO_i)\}
\end{equation}

\end{lem}

\begin{proof}

\textit{Definition of $u^\diamond$}: Define $J(v)\colonequals\{i\in \mathbb{S}:\langle v,\varpi_{\calO_i}\rangle\not= \height_{L_0}(\calO_i)\}$ and let $K$ be the connected component of $J(v)$ containing $k$. Applying \Cref{subdiag} to $K$ and $k$, we call the resulting $\sum_{i\in K}c_i\alpha_i^\vee$ part $u$ and denote its $\sigma$-average as $u^\diamond$. We want to set $v'=v-\epsilon u^\diamond$ for some $\epsilon\in\mathbb{Q}_{>0}$.

\textit{Construction of $v'$}: Note that $\langle v,\varpi_{\calO_i}\rangle -\height_{L_0}(\calO_i)>0$ for all $i\in K$ by definition of $J(v)$. As $\langle u^\diamond,\varpi_{\calO_i}\rangle=\sum_{j\in\calO_i}c_j>0$ for all $i\in K$ by \Cref{subdiag}, we can find the maximal one $\epsilon_1$ such that $\langle v-\epsilon_1 u^\diamond,\varpi_{\calO_i}\rangle-\height_{L_0}(\calO_i)\ge0$ for all $i\in K$. We set $\epsilon\colonequals\max\{\epsilon_1,\#\calO_{k}\cdot \langle v,\alpha_{k}\rangle \}$ and $v'\colonequals v-\epsilon u^\diamond$.

\textit{Verification of $v'\in\renv{L_0}$}: Observing that $\langle u^\diamond,\alpha_i\rangle>0$ only when $i\in\calO_{k}$ and $\langle u^\diamond,\alpha_i\rangle<0$ only when $i\in \partial K$, we only need to consider $k' \in\calO_{k}$. Since $v$ is $\sigma$-invariant, we have $\langle v-\epsilon u^\diamond,\alpha_{k'}\rangle=\langle v,\alpha_{k}\rangle-\frac{1}{\#\calO_{k}}\epsilon$, which is non-negative by the definition of $\epsilon$. As $u^\diamond$ is $\sigma$-invariant, so is $v'$.

\textit{Verification of ``\eqref{inS} or \eqref{breakpts}''}: As $K\subset J(v)$ and $\langle u^\diamond,\varpi_{\calO_i}\rangle=0~\forall i\in \mathbb{S}\setminus\! K$,\begin{align*}\{i\in\mathbb{S}:\langle v,\varpi_{\calO_i}\rangle=\height_{L_0}(\calO_i)\}&=\{i\in \mathbb{S}\setminus\! K:\langle v,\varpi_{\calO_i} \rangle =\height_{L_0}(\calO_i)\}\\&=\{i\in\mathbb{S}\setminus\! K:\langle v',\varpi_{\calO_i}\rangle=\height_{L_0}(\calO_i)\}.
\end{align*} So, \eqref{inS} is equivalent to $\{i\in K:\langle v-\epsilon u^\diamond,\varpi_{\calO_i}\rangle=\height_{L_0}(\calO_i)\}\neq\emptyset$ which is equivalent to that $\epsilon=\epsilon_1$. Next, from the previous paragraph, we know that $\mathbb{S}\setminus\! I(v')=(\mathbb{S}\setminus\! I(v) \setminus\!\calO_{k})\cup \partial K$ if and only if $\epsilon=\#\calO_{k}\cdot\langle v,\alpha_{k}\rangle$.\footnote{If $\epsilon_1=\epsilon<\#\calO_{k}\cdot\langle v,\alpha_{k}\rangle$, then $\mathbb{S}\setminus\! I(v')=(\mathbb{S}\setminus\! I(v))\cup\partial K$.} However, for $i\in\partial K$, we have $\langle v,\varpi_{\calO_i}\rangle=\langle v',\varpi_{\calO_i}\rangle=\height_{L_0}(\calO_i)$ but $k$ belongs to only the right-hand side set in \eqref{breakpts}. This proves the claim.
\end{proof}
\begin{lem}\label{vcannotbezero}
    Let $v\in B(G,\mu)$. Then $\langle v,\varpi_{\calO_i}\rangle =0$ implies that $i\in I(v)$.
\end{lem}
\begin{proof}
Suppose that $\langle v,\varpi_{\calO_i}\rangle = 0$. Since $v$ is $\sigma$-invariant, we have $\langle v,\varpi_i\rangle = 0$. Note that $\varpi_i \in \sum_{j\in\mathbb{S}} \mathbb{Q}_{\ge0}\alpha_j$ and $\langle v,\alpha_i\rangle\ge0$ as $v$ is dominant. Hence $\langle v,\alpha_i\rangle =0$.
\end{proof}

\begin{proof}[Proof of \Cref{equivsmall}]

    1)$\Rightarrow$2): Trivial.
    
    2)$\Rightarrow$3): Assume that 3) does not hold. As $v\in B(G,\mu)_{\ge L_0}\subset\renv{L_0}$, we can apply \Cref{goingdown}. If the new $v'$ satisfies the assumption again, we keep  repeating this process. As the number in \eqref{inS} is $\le\# S$ and that in \eqref{breakpts} is $\ge0$, the process terminates in finite steps. The final $v_{\fin}$ then satisfies that $\langle v_{\fin},\varpi_{\calO_i}\rangle= \height_{L_0}(\calO_i)$ for all $i\in \mathbb{S}\setminus\! I(v_{\fin})$. By \Cref{vcannotbezero}, $\height_{L_0}(\calO_i)\neq0$ and hence $\height_{L_0}(\calO_i)=\langle \mu,\varpi_{\calO_i}\rangle -m_i$ for some $m_i\in\mathbb{Z}_{\ge1}$. We then have $\langle\mu- v_{\fin},\varpi_{\calO_i}\rangle\in\mathbb{Z}_{\ge1}$ for all $i\in \mathbb{S}\setminus\! I(v_{\fin})$, which means that $v_{\fin}\in B(G,\mu)_{\ge L_0}$ and so $v$ was not minimal. It contradicts to 2).

    3)$\Rightarrow$1): Suppose that $w\in B(G,\mu)_{\ge L_0}$. Then $\langle w, \varpi_i\rangle \ge \height_{L_0}(\calO_i) = \langle v, \varpi_{\calO_i}\rangle$ for all $i \in \mathbb{S}\setminus\!I(v)$ by the assumption. We get $\langle w-v,\varpi_{i}\rangle\ge0$ for all $i\in\mathbb{S}\setminus\! I(v)$ since $w$ and $v$ are $\sigma$-invariant.
    
    Let $j\in I(v)$. The dual version of \Cref{subdiag} applied to the connected component of $I(v)$ containing $j$ tells us that $\varpi_j$ is a nonnegative linear combination of $\alpha_i$ ($i\in I(v)$) and $\varpi_k$ ($k\in \mathbb{S}\setminus\! I(v)$). However, $\langle w,\alpha_i\rangle\ge0$ as $w\in V^+$ and $\langle v,\alpha_i\rangle=0$ for all $i\in I(v)$, so we have $\langle w-v,\alpha_i\rangle\ge0$ for all $i\in I(v)$. Therefore, $\langle w-v,\varpi_j\rangle\ge0$ for all $j\in I(v)$.\\

Now we prove that the smallest element $v\colonequals B(G,\mu)_{\ge L_0}$ lies in $\bgmuindec$. Suppose $v\not\in \bgmuindec$. Then there exists a $\sigma$-stable $J\subsetneq\mathbb{S}$ such that $I(v) \subset J$ and $\mu - v \in \sum_{j\in J}\mathbb{R}_{\ge0} \alpha_j^\vee$. Hence, we get $\langle \mu - v, \varpi_{\calO_i}\rangle = 0$ for all $i \not\in J$. As $J$ is a proper subset of $\mathbb{S}$, we can choose $i_0\in \mathbb{S}\setminus\! J\subset \mathbb{S}\setminus\!I(v)$.
By (3), this implies that $\langle v,\varpi_{\calO_{i_0}}\rangle = \height_{L_0}(\calO_{i_0})$. However, $\height_{L_0}(\calO_i) <  \langle\mu,\varpi_{\calO_{i}}\rangle$ for any $i\in\mathbb{S}$ by the construction of $L$ and the height function. This gives a contradiction and so $v\in \bgmuindec$.
\end{proof}

\section{Lemmas for \Cref{motiv}}\label{cutelemmas}

\subsection{Notations}
A polygon is assumed to have \textit{lattice points as its vertices}. Denote by $\area(P)$ the area of $P$, by $\iota(P)$ and $\beta(P)$ the number of lattice points interior to $P$ and on the boundary of $P$. Denote $O\colonequals (0,0)$, $X\colonequals(i,0)$, and $Y\colonequals(i,n-i)$. Denote $\tri$ as the triangle $\triangle OXY$.
\subsection{Two lemmas}

\begin{lem}\label{btob-a}
    Let $\calC_k$ and $\calD_k$ be the sets defined in \Cref{convexhull}. Then
    \[\sum_{\substack{k\ge1\\ ((a_l,b_l))_{l\le k}\in \square_k}} (\q-1)^{k-1}\q^{1-k+\frac{\sum_{1\le l_1<l_2\le k}(a_{l_1}b_{l_2}-a_{l_2}b_{l_1}) +\sum_{1\le l\le k}\gcd(a_l,b_l)}{2}}\] for $\square=\calC$ and $\square =\calD$ are equal. Moreover, $\calC_k$ can be identified with the set of convex broken lines connecting $O$ and $Y$ with $k$ segments and lying strictly above $\overline{OXY}$.
\end{lem}

\begin{proof}
    The first part follows easily under the map $\calC_k\to\calD_k$ given by $(x_l,y_l)\mapsto (x_l,x_l+y_l)$, and the second part is straightforward using the correspondence $\calC_k\ni ((x_l,y_l))_{l\le k}\mapsto$ the broken line connecting $(0,0)$, $(x_1,y_1)$, $(x_1+x_2,y_1+y_2)$, $\cdots$, $(x_1+\cdots+x_k, y_1+\cdots+y_k)$.
\end{proof}

\begin{lem}\label{pick} Regarding $\calC_k$ as the set of broken lines from \Cref{btob-a}, \Cref{motiv} is equivalent to\[\sum_{k\ge1,~C\in\calC_k}\left(\frac{\q-1}{\q}\right)^{k-1}\q^{-u(C)}=1,\]where $u(C)$ is the number of lattice points below $C$ and above the $x$-axis.
\end{lem}

\begin{proof}
For $C\in\calC_k$, let $P_C \colonequals C\cup \overline{OY}$ the polygon corresponding to $C$.

For the left-hand side of \Cref{motiv}, it is easy to show by induction that
\[\frac{1}{2}\sum_{1\le l_1<l_2\le k}(x_{l_1}y_{l_2}-x_{l_2}y_{l_1})=\area(P_C).\]The number of lattice points on the line connecting $(m,n)\in\mathbb{Z}^2$ and $(m+a,n+b)$ is $\gcd(a,b)+1$;
\[\sum_{1\le l\le k}\gcd(x_l,y_l)=\beta(P_C)-\gcd(i,n-i).\]

For the right-hand side, we know $\area(\tri)=\frac{i(n-i)}{2}$ and $\beta(\tri)=n+\gcd(i,n-i)$. Now we apply the well-known Pick's Theorem: \textit{For every (lattice point) polygon $P$, we have $A(P) = \iota(P) + \frac{\beta(P)}{2} -1$.}

Then we get\[\sum_{k\ge1,~ C\in \calC_k} (\q-1)^{k-1}\q^{-(k-1)+\iota(P_C)+\beta(P_C)}=\q^{\iota(\tri)+\beta(\tri)-(n-1)}.\]

Observe that, for $C\in\calC_k$, we have $u(C) = \iota(\tri) +\beta(\tri) -(n-1) - \iota(P_C) -\beta(P_C)$ and $b(C) = k-1$. Since $\bigcup_k\calC_k = C$, we get the conclusion.
\end{proof}

\begin{misc}
    
\begin{rem}
We do not know if the identity such as \Cref{essentialidea} is well-known. It looks interesting to us because the left-hand side is not homogeneous in the sense that $u(C)+b(C)$ is not constant but it gives a way to generate $1$ using polynomials of the form $\p^a(1-\p)^b$. We wonder if the set $\{(u(C),b(C)):C\in \calC\}$ parametrizes all such pairs $\{(a_i,b_i)_i\}$.

More precisely, let $S$ be a finite set (possibly with dulicates) of $(a,b)$ such that $a$ and $b$ are natural numbers and containing $(0,1)$. Suppose that\[\sum_{(a,b)\in S} \p^a (1-\p)^b = 1.\]

Then we would like to ask if there exist $m,n\in\mathbb{N}$ such that $S=\{(u(P),v(P)-2):P\in C_{m,n}\}$ where $C_{m,n}$ is the set defined in \Cref{mainlem} corresponding to the triangle $\tri_{m,n}$ whose vertices are $(0,0)$, $(m,0)$, and $(m,n)$.
\end{rem}
\end{misc}

\begin{misc}
\begin{prop}[The motivation]\label{motiv}
Fix natural numbers $i< n$, then
\[\sum_{\substack{k\ge1\\ ((a_l,b_l))_{l\le k}\in D_k}} (\mathbf{q}-1)^{k-1}\mathbf{q}^{1-k+\frac{\sum_{1\le l_1<l_2\le k}(a_{l_1}\!b_{l_2}-a_{l_2}\!b_{l_1}) +\sum_{1\le l\le k}\gcd(a_l,b_l)}{2}}=\mathbf{q}^{\frac{i(n-i)-n}{2}+1},\]where $D_k=\{((a_l,b_l))_{l\le k}: a_l,b_l\text{ are natural numbers and } a_1+\cdots+ a_k=i,~b_1+\cdots+b_k=n,\text{ and } 1>\frac{a_1}{b_1}>\cdots>\frac{a_k}{b_k}>0.\}$.
\end{prop}
\end{misc}

\begin{misc}
In fact, the above set is called HN-indecomposable is closely related to HN-irreducible in the following sense. When the root system has a single $\sigma$-orbit of connected components,  any HN-indecomposable pair $(\mu,v)$ is either HN-irreducible, or satisfies $v=\mu$ with $\mu$ central (see \cite[Theorem 2.5.6]{CKV}). More generally, for each $\sigma$-orbit of connected components, an HN-indecomposable pair $(\mu,v)$ is either HN-irreducible or has $v=\mu$ central.

However, f
\end{misc}

\bibliographystyle{alpha}
\bibliography{hnybib}

\end{document}